# Recurrence of cocycles and stationary random walks

## Klaus Schmidt[1,2]

*University of Vienna and Erwin Schrödinger Institute*

**Abstract:** We survey distributional properties of $\mathbb{R}^d$-valued cocycles of finite measure preserving ergodic transformations (or, equivalently, of stationary random walks in $\mathbb{R}^d$) which determine recurrence or transience.

Let $(X_n, n \geq 0)$ be an ergodic stationary $\mathbb{R}^d$-valued stochastic process, and let $(Y_n = X_0 + \cdots + X_{n-1}, n \geq 1)$ be the associated random walk. What can one say about recurrence of this random walk if one only knows the distributions of the random variables $Y_n$, $n \geq 1$?

It turns out that methods from ergodic theory yield some general sufficient conditions for recurrence of such random walks without any assumptions on independence properties or moments of the process $(X_n)$.

Most of the results described in this note have been published elsewhere. Only the Theorems 12 and 14 on recurrence of symmetrized random walks are — to my knowledge — new.

Let us start our discussion by formulating the recurrence problem in the language of ergodic theory.

Let $T$ be a measure preserving ergodic automorphism of a probability space $(X, \mathsf{S}, \mu)$, $d \geq 1$, and let $f\colon X \longrightarrow \mathbb{R}^d$ be a Borel map. The *cocycle* $\mathbf{f}\colon \mathbb{Z} \times X \longrightarrow \mathbb{R}^d$ is given by

$$\mathbf{f}(n,x) = \begin{cases} f(T^{n-1}x) + \cdots + f(x) & \text{if } n > 0, \\ 0 & \text{if } n = 0, \\ -\mathbf{f}(-n, T^n x) & \text{if } n < 0, \end{cases} \qquad (1)$$

and satisfies that

$$\mathbf{f}(m, T^n x) + \mathbf{f}(n, x) = \mathbf{f}(m+n, x)$$

for all $m, n \in \mathbb{Z}$ and $x \in X$.

**Definition 1.** The cocycle $\mathbf{f}$ (or the function $f$) in (1) is *recurrent* if

$$\liminf_{n \to \infty} \|\mathbf{f}(n,x)\| = 0$$

for $\mu$-a.e. $x \in X$, where $\|\cdot\|$ is a norm on $\mathbb{R}^d$.

If we start with a stationary $\mathbb{R}^n$-valued random walk $(X_n)$ on a probability space $(\Omega, \mathcal{T}, P)$ we may assume without loss in generality that the sequence of random variables $X_n\colon \Omega \longrightarrow \mathbb{R}^d$ is two-sided and generates the sigma-algebra $\mathcal{T}$.

---

[1] Mathematics Institute, University of Vienna, Nordbergstraße 15, A-1090 Vienna, Austria.
[2] Erwin Schrödinger Institute, for Mathematical Physics, Boltzmanngasse 9, A-1090 Vienna, Austria, e-mail: klaus.schmidt@univie.ac.at

Then there exists an ergodic, probability-preserving automorphism $T$ of $(\Omega, \mathcal{T}, P)$ with $X_{m+n} = X_m \circ T^n$ for all $m, n \in \mathbb{Z}$, and we set $f = X_0$ and obtain that $\mathbf{f}(n, \cdot) = Y_n = X_0 + \cdots + X_{n-1}$ for every $n \geq 1$. Our notion of recurrence for $f$ coincides with the usual probabilistic notion of recurrence of the random walk $(Y_n, n \geq 1)$.

We return to the ergodic-theoretic setting. The cocycle $\mathbf{f}$ (or the function $f$) is a *coboundary* if there exists a Borel map $b \colon X \longrightarrow \mathbb{R}^d$ such that

$$f(x) = b(Tx) - b(x) \text{ for } \mu\text{-a.e. } x \in X \tag{2}$$

or, equivalently, if

$$\mathbf{f}(n, x) = b(T^n x) - b(x)$$

for every $n \in \mathbb{Z}$ and $\mu$-a.e. $x \in X$.

If $f' \colon X \longrightarrow \mathbb{R}^d$ is a second Borel map and $\mathbf{f}' \colon \mathbb{Z} \times X \longrightarrow \mathbb{R}^d$ the resulting cocycle in (1), then $\mathbf{f}$ and $\mathbf{f}'$ (or $f$ and $f'$) are *cohomologous* if there exists a Borel map $b \colon X \longrightarrow \mathbb{R}^d$ with

$$f'(x) = b(Tx) + f(x) - b(x) \text{ for } \mu\text{-a.e. } x \in X. \tag{3}$$

Recurrence is easily seen to be a *cohomology invariant*.

**Proposition 2.** *Let $T$ be an ergodic measure preserving automorphism of a standard probability space $(X, \mathcal{S}, \mu)$, and let $f, f' \colon X \longrightarrow \mathbb{R}^d$ be Borel maps. Then $f$ is recurrent if and only if $f'$ is recurrent.*

The proof of Proposition 2 can be found, for example, in [8] (cf. also [9]). The key to this proposition is the observation that recurrence is equivalent to the following condition: for every $B \in \mathcal{S}$ with $\mu(B) > 0$ and every $\varepsilon > 0$ there exists a nonzero $n \in \mathbb{Z}$ with

$$\mu(B \cap T^{-n} B \cap \{x \in X : T^n x \neq x \text{ and } \|\mathbf{f}(n, x)\| < \varepsilon\}) > 0. \tag{4}$$

The recurrence properties of $\mathbb{R}^n$-valued random walks arising from i.i.d. processes are understood completely. For random walks arising from ergodic stationary processes the question of recurrence is much more complex. The first results in this direction appear to be due to [1] and [6].

**Theorem 3.** *Let $T$ be an ergodic measure preserving automorphism of a standard probability space $(X, \mathcal{S}, \mu)$ and $f \colon X \longrightarrow \mathbb{R}$ a Borel map. If $f$ is integrable, then it is recurrent if and only if $\int f \, d\mu = 0$.*

Integrable functions with zero integral satisfy the strong law of large numbers by the ergodic theorem. For real-valued functions satisfying the weak law of large numbers, recurrence was proved in [3], [9] and by B. Weiss (unpublished).

**Theorem 4.** *Let $T$ be an ergodic measure preserving automorphism of a standard probability space $(X, \mathcal{S}, \mu)$. If a Borel map $f \colon X \longrightarrow \mathbb{R}$ satisfies the weak law of large numbers (i.e. if $\lim_{n \to \infty} \mathbf{f}(n, \cdot)/n = 0$ in measure), then $f$ is recurrent.*

The failure of $f$ to be recurrent may be due to very simple reasons. For example, if $d = 1$ and $f$ is integrable with nonzero integral, then $f$ is nonrecurrent, but by subtracting the integral we add a *drift term* to the cocycle $\mathbf{f}$ which makes it recurrent. This motivates the following definition.

**Definition 5.** The *recurrence set* of a Borel map $f\colon X \longrightarrow \mathbb{R}^d$ is defined as

$$R(f) = \{c \in \mathbb{R}^d : f - c \text{ is recurrent}\}.$$

It is not difficult to see that $R(f)$ is a Borel set (cf. [9]). In order to discuss $R(f)$ further we consider the *skew-product transformation* $T_f\colon X \times \mathbb{R}^d \longrightarrow X \times \mathbb{R}^d$ defined by

$$T_f(x,g) = (Tx, f(x) + g) \tag{5}$$

for every $(x,g) \in X \times \mathbb{R}^d$. It is clear that $T_f$ preserves the infinite measure $\mu \times \lambda$, where $\lambda$ is the Lebesgue measure on $\mathbb{R}^d$.

For the following result we recall that a Borel set $D \subset X \times \mathbb{R}^d$ is *wandering* under $T_f$ if $T_f^n D \cap D = \varnothing$ for every nonzero $n \in \mathbb{Z}$. The transformation $T_f$ is *conservative* if every $T_f$-wandering set $D \subset X \times \mathbb{R}^d$ satisfies that $(\mu \times \lambda)(D) = 0$.

**Proposition 6 ([9]).** *Let $T$ be an ergodic measure preserving automorphism of a standard probability space $(X, \mathcal{S}, \mu)$ and $f\colon X \longrightarrow \mathbb{R}^d$ a Borel map. Then $f$ is recurrent if and only if $T_f$ is conservative. If $f$ is transient (i.e. not recurrent), then there exists a $T_f$-wandering set $D \subset X \times \mathbb{R}^d$ with*

$$(\mu \times \lambda)\left((X \times \mathbb{R}^d) \smallsetminus \bigcup_{n \in \mathbb{Z}} T_f^n D\right) = 0 \tag{6}$$

*(such a set is sometimes called a* sweep-out set*)*.

Since recurrence is a cohomology invariant, we obtain the following result.

**Corollary 7.** *The recurrence set is a cohomology invariant. In particular, $R(f) = R(f \circ T)$.*

Here are some examples of recurrence sets.

**Examples 8.** Let $f\colon X \longrightarrow \mathbb{R}$ be a Borel map.
 (1) If $f$ is integrable then $R(f) = \{\int f \, d\mu\}$.
 (2) If $f \geq 0$ and $\int f \, d\mu = \infty$ then $R(f) = \varnothing$.
 (3) If $(X_n)$ is the i.i.d. Cauchy random walk and $f = X_0$ then $R(f) = \mathbb{R}$.
 (4) The recurrence set $R(f)$ can be equal to any given countable closed subset of $\mathbb{R}$ (unpublished result by B. Weiss).
 (5) The sets $R(f)$ and $\mathbb{R} \smallsetminus R(f)$ can simultaneously be dense in $\mathbb{R}$. An example appears in [9]: let $T$ be the tri-adic adding machine on $X = \{0,1,2\}^{\mathbb{N}}$, and let $\phi\colon X \longrightarrow X$ be the map which interchanges the digits '1' and '2' in each coordinate. Then $\phi^2 =$ Identity and the automorphisms $T$ and $T' = \phi \circ T \circ \phi$ have the same orbits in the complement of a $\mu$-null set ($\mu$ is the Haar measure of $X$). Hence there exists a Borel map $f\colon X \longrightarrow \mathbb{Z}$ with $Tx = T'^{f(x)}x$ for $\mu$-a.e. $x \in X$, and both $R(f)$ and $\mathbb{R} \smallsetminus R(f)$ are dense in $\mathbb{R}$.

Corollary 7 and Example 8 (1) shows that the recurrence set behaves somewhat like an invariant integral for possibly nonintegrable functions (except for linearity). Furthermore, if $R(f) = \varnothing$ or $|R(f)| > 1$ (where $|\cdot|$ denotes cardinality), then $f$ cannot be cohomologous to an $L^1$-function.

In [9] the problem was raised how one could recognize whether the recurrence set of a function $f$ is nonempty. For real-valued functions a partial answer to this question was given in [5].

**Theorem 9.** *Let $T$ be an ergodic measure preserving automorphism of a standard probability space $(X, \mathcal{S}, \mu)$ and $f\colon X \longrightarrow \mathbb{R}$ a Borel map. If there exist $\varepsilon, K > 0$ such that $\mu(\{x \in X : |\mathbf{f}(n,x)| \leq K\}) > \varepsilon$ for every $n \geq 1$, then $R(f) \neq \varnothing$.*

*In particular, if the distributions of the random variables $\mathbf{f}(n,\cdot)$, $n \geq 1$, are uniformly tight, then $R(f) \neq \varnothing$.*

The proof of Theorem 9, as well as of several related results, depends on a 'local limit formula' in [10]:

Let $f$ take values in $\mathbb{R}^d$, $d \geq 1$, and let $\sigma_k^{(d)}$ be the distribution of the map $\mathbf{f}(k,\cdot)/k^{1/d}$ and $\tau_k^{(d)} = \frac{1}{k}\sum_{l=1}^k \sigma_l^{(d)}$. If $f$ is transient there exist an integer $L \geq 1$ and an $\varepsilon > 0$ such that

$$\limsup_{k\to\infty} \tau_k^{(d)}(B(\eta)) \leq 2^d L \varepsilon^{-d} \lambda(B(\eta)) \tag{7}$$

and

$$\limsup_{k\to\infty} \sum_{n=0}^N 2^n \tau_{2^n k}^{(d)}(B(2^{-n/d}\eta)) \leq 2^{d+1} d L^d \varepsilon^{-d} \lambda(B(\eta)) \tag{8}$$

for every $\eta > 0$ and $N \geq 1$, where $\|\cdot\|$ denotes the maximum norm on $\mathbb{R}^d$ and

$$B(\eta) = \{v \in \mathbb{R}^d : \|v\| < \eta.\}$$

The proof of these formulae uses abstract ergodic theory (existence of a sweep-out set for the skew-product, cohomology and orbit equivalence).

The inequality (7) shows that the possible limits of the sequences $(\tau_k^{(d)}, k \geq 1)$ are absolutely continuous at 0 and gives a bound on their density functions there. As a corollary we obtain that $f$ must be recurrent if any limit point of $(\tau_k^{(d)}, k \geq 1)$ has an atom at 0 (which proves Theorem 4, for example).

In order to give a very scanty idea of the proof of Theorem 9 we choose an increasing sequence $(k_m)$ such that the vague limits $\rho_n = \lim_{m\to\infty} \tau_{2^n k_m}^{(d)}$, $n \geq 1$, exist, and obtain from (8) that

$$\sum_{n=0}^N 2^n \rho_n(B(2^{-n/d}\eta)) \leq 2^{d+1} d L^d \varepsilon^{-d} \lambda(B(\eta))$$

for every $\eta > 0$. This shows that some of the $\rho_n$ must have arbitrarily small density at 0.

If, under the hypotheses of Theorem 4, $R(f)$ were empty, one could construct limit points of certain averages of the $\tau_k^{(1)}$ with arbitrarily small total mass, which would violate the hypotheses of Theorem 9. The details can be found in [5].

For $d = 2$, these considerations lead to the following special case of a result in [10].

**Theorem 10 ([2]).** *Let $T$ be an ergodic measure preserving automorphism of a standard probability space $(X, \mathcal{S}, \mu)$ and $f\colon X \longrightarrow \mathbb{R}^2$ a Borel map. If the distributions of the functions $\mathbf{f}(n,\cdot)/\sqrt{n}$, $n \geq 1$, in (1) converge to a Gaussian distribution on $\mathbb{R}^2$, then $f$ is recurrent.*

One might be tempted to conjecture that, for a Borel map $f\colon X \longrightarrow \mathbb{R}^2$, uniform tightness of the distributions of the random variables $(\mathbf{f}(n,\cdot)/\sqrt{n}, n \geq 1)$, would imply recurrence. However, an example by M. Dekking in [3] shows that the distributions of the functions $\mathbf{f}(n,\cdot)$, $n \geq 1$, can be uniformly tight even if $f$ is transient.

Let me mention a version of Theorem 10 for $\mathbb{R}^d$-valued maps.

**Theorem 11.** *Let $T$ be an ergodic measure preserving automorphism of a standard probability space $(X, \mathcal{S}, \mu)$ and $f\colon X \longrightarrow \mathbb{R}^d$ a Borel map, where $d \geq 1$. Suppose that there exists an $\varepsilon > 0$ such that every limit point $\rho$ of the distributions of the random variables $\mathbf{f}(n, \cdot)/n^{1/d}$, $n \geq 1$ satisfies that*

$$\liminf_{\eta \to 0} \rho(B(\eta))/\eta^d \geq \varepsilon.$$

*Then $f$ is recurrent.*

This result may not look very interesting for $d > 2$, but it (or at least an analogue of it) can be used to prove recurrence and ergodicity of skew-product extensions for cocycles with values in certain (noncommutative) matrix groups, such as the group of unipotent upper triangular $d \times d$ matrices (cf. [4]).

Theorem 10 can be used to prove *recurrence* of an $\mathbb{R}^2$-valued function, but I don't know of any useful information about *nonemptiness of the recurrence set* of a transient function $f\colon X \longrightarrow \mathbb{R}^2$. A possible approach to this problem is to investigate recurrence of the 'symmetrized' version $\tilde{f}\colon X \times X \longrightarrow \mathbb{R}^d$ of $f$, defined by

$$\tilde{f}(x, y) = f(x) - f(y) \tag{9}$$

for every $(x, y) \in X \times X$. We denote by $\tilde{\mathbf{f}}\colon \mathbb{Z} \times (X \times X) \longrightarrow \mathbb{R}^d$ the cocycle for the transformation $S = T \times T$ on $X \times X$ defined as in (1) (with $S$ and $\tilde{f}$ replacing $T$ and $f$).

For $d = 1$ we have the following result.

**Theorem 12.** *Let $T$ be an ergodic measure preserving automorphism of a standard probability space $(X, \mathcal{S}, \mu)$ and $f\colon X \longrightarrow \mathbb{R}$ a Borel map. Suppose that the distributions of the random variables $\mathbf{f}(n, \cdot)/n$, $n \geq 1$, are uniformly tight. Then the map $\tilde{f}\colon X \times X \longrightarrow \mathbb{R}$ in (9) is recurrent.*

*Conversely, if the distributions of the random variables $\tilde{\mathbf{f}}(n, \cdot)/n$, $n \geq 1$, are uniformly tight, then $R(f) \neq \varnothing$.*

*Proof.* Let $\sigma_n = \sigma_n^{(1)}$ be the distribution of the random variable $\mathbf{f}(n, \cdot)/n, n \geq 1$, and let $h_\delta = 1_{[-\delta/2, \delta/2]}\colon \mathbb{R} \longrightarrow \mathbb{R}$ be the indicator function of the interval $[-\delta/2, \delta/2] \subset \mathbb{R}$ for every $\delta \in (0, 1)$. We set $g_\delta = \frac{1}{\delta^2} \cdot h_\delta * h_\delta$, where $*$ denotes convolution. Then $\int g_\delta \, d\lambda = 1$, where $\lambda$ is Lebesgue measure on $\mathbb{R}$,

$$0 \leq g_\delta \leq \tfrac{1}{\delta} \cdot 1_{[-\delta, \delta]} \tag{10}$$

and

$$1 = \iiint g_\delta(x + y + z) \, d\sigma_n(x) \, d\sigma_n(-y) \, dz$$
$$= \frac{1}{\delta^2} \cdot \iint (h_\delta * \sigma_n)(u + z)(h_\delta * \sigma_n)(u) \, du \, dz.$$

We put

$$\phi_\delta(z) = \iint g_\delta(x + y + z) \, d\sigma_n(x) \, d\sigma_n(-y)$$

for every $z \in \mathbb{R}$.

Assume that the probability measures $(\sigma_n, n \geq 1)$ are uniformly tight and choose $K > 0$ so that $\sigma_n([-K/2, K/2]) > 1/2$ for all $n \geq 1$. Then

$$\int_{-K-1}^{K+1} \phi_\delta(u) \, du \geq (\sigma_n * \sigma_n)([-K, K]) > 1/4$$

for every $n \geq 1$. Hence

$$\lambda(\{u \in \mathbb{R} : \phi_\delta(u) > 1/(8K+8)\}) > 0.$$

Since $\phi_\delta$ assumes its maximum at 0 we conclude that

$$\phi_\delta(0) > 1/(8K+8) \tag{11}$$

for every $\delta > 0$.

We set $\bar{\sigma}_n(B) = \sigma_n(-B)$ for every Borel set $B \subset \mathbb{R}$, denote by $\tilde{\sigma}_n = \sigma_n * \bar{\sigma}_n$ the distribution of the random variable $\tilde{\mathbf{f}}(n, \cdot)/n$, and obtain that

$$2^k \cdot \tilde{\sigma}_n([-2^{-k}\eta, 2^{-k}\eta]) \geq \eta \cdot \int g_{2^{-k}\eta} d\tilde{\sigma}_n \geq \eta/(8K+8)$$

for every $\eta > 0$ and $k, n \geq 1$, by (10)–(11). By comparing this with (8) we see that $\tilde{f}$ is recurrent.

Now suppose that the distributions $\tilde{\sigma}_n$ of the $\tilde{\mathbf{f}}(n, \cdot)/n$, $n \geq 1$, are uniformly tight. Theorem 2.2 in [7] implies the existence of a sequence $(a_n, n \geq 1)$ in $\mathbb{R}$ such that the probability measures $p_{a_n} * \sigma_n$, $n \geq 1$, are uniformly tight, where $p_t$ is the unit point-mass at $t$ for every $t \in \mathbb{R}$.

We set $\sigma_0 = p_0$, $a_0 = 0$, and $\sigma_{-n} = \bar{\sigma}_n$ and $a_{-n} = -a_n$ for every $n \geq 1$, and obtain that the family of probability measures $\{p_{a_n} * \sigma_n : n \in \mathbb{Z}\}$ is uniformly tight. It follows that the map $(m, n) \mapsto c(m, n) = a_{m+n} - a_m - a_n$ from $\mathbb{Z} \times \mathbb{Z}$ to $\mathbb{R}$ is bounded, and we choose a translation-invariant mean $M(\cdot)$ on $\ell^\infty(\mathbb{Z}, \mathbb{R})$ and set $b_n = M(c(n, \cdot))$ for every $n \in \mathbb{Z}$. Then the set $\{b_n : n \in \mathbb{Z}\}$ is bounded,

$$c(m, n) = b_{m+n} - b_m - b_n$$

for every $m, n \in \mathbb{Z}$, and the maps $n \mapsto a_n$ and $n \mapsto b_n$ differ by a homomorphism from $\mathbb{Z}$ to $\mathbb{R}$ of the form $n \mapsto tn$ for some $t \in \mathbb{R}$.

Our choice of $t$ implies that the distributions of the random variables $\mathbf{f}(n, \cdot)/n + a_n - b_n = \mathbf{f}(n, \cdot)/n + tn$, $n \geq 1$, are uniformly tight. For every $\varepsilon > 0$ we can therefore choose a $K > 0$ such that

$$\mu(\{x \in X : -nK - tn^2 \leq \mathbf{f}(n, x) \leq nK - tn^2\}) > 1 - \varepsilon$$

for every $n \in \mathbb{Z}$. If $t \neq 0$ we obtain a contradiction for sufficiently small $\varepsilon$. This proves that the sequence $(\rho_n, n \geq 1)$ is uniformly tight and that $R(f) \neq \varnothing$ by Theorem 9. □

**Remarks 13.** (1) Under the hypotheses of Theorem 12, $R(\tilde{f}) \neq \varnothing$ by Theorem 9. It should really be obvious that in this case $0 \in R(\tilde{f})$, but I can't think of any direct reason why this should be true.

(2) Let $T$ be an ergodic measure preserving automorphism of a standard probability space $(X, \mathcal{S}, \mu)$ and $f : X \longrightarrow \mathbb{R}$ a Borel map. If the distributions of the random variables $(\mathbf{f}(n, \cdot), n \geq 1)$ are uniformly tight, then Theorem 12 implies that there exists, for $(\mu \times \mu)$-a.e. $(x, y) \in X \times X$, an increasing sequence $(n_k, k \geq 1)$ of natural numbers with $\lim_{k \to \infty} \mathbf{f}(n_k, x) - \mathbf{f}(n_k, y) = 0$.

(3) From the proof of Theorem 12 it is clear that we can replace the uniform tightness of the distributions of the $\mathbf{f}(n, \cdot)$, $n \geq 1$, by the condition that there exist $\varepsilon, K > 0$ such that $\mu(\{x \in X : |\mathbf{f}(n, x)| \leq K\} > \varepsilon)$ for every $n \geq 1$.

We turn to the much more interesting case where $d > 1$. For a function $f : X \longrightarrow$

$\mathbb{R}^d$, Dekking's example in [3] shows that uniform tightness of the distributions of of the random variables $\mathbf{f}(n,\cdot)/n^{1/d}$, $n \geq 1$, need not imply recurrence. However, the function $\tilde{f}$ in (9) *is* recurrent under this condition. At present there appears to be no analogue of the reverse implication of Theorem 12.

**Theorem 14.** *Let $T$ be an ergodic measure preserving automorphism of a standard probability space $(X, \mathcal{S}, \mu)$ and $f \colon X \longrightarrow \mathbb{R}^d$ a Borel map. Suppose that the distributions of the random variables $\mathbf{f}(n,\cdot)/n^{1/d}$, $n \geq 1$, are uniformly tight (or that there exist $\varepsilon, K > 0$ such that $\mu(\{x \in X : \|\mathbf{f}(n,x)/n^{1/d}\| \leq K\}) > \varepsilon$ for every $n \geq 1$). Then the map $\tilde{f} \colon X \times X \longrightarrow \mathbb{R}^d$ in (9) is recurrent.*

*Proof.* The proof is essentially the same as that of the first part of Theorem 12. Let $\sigma_n = \sigma_n^{(d)}$ be the distribution of the random variable $\mathbf{f}(n,\cdot)/n^{1/d}$, and let $h_\delta = 1_{[-\delta/2,\delta/2]^d} \colon \mathbb{R}^d \longrightarrow \mathbb{R}$ be the indicator function of $[-\delta/2, \delta/2]^d \subset \mathbb{R}^d$ for every $\delta \in (0,1)$. We set $g_\delta = \frac{1}{\delta^{2d}} \cdot h_\delta * h_\delta$ and put

$$\phi_\delta(z) = \iint g_\delta(x + y + z)\, d\sigma_n(x)\, d\sigma_n(-y)$$

for every $z \in \mathbb{R}^d$. As in the proof of Theorem 12 we see that

$$2^{dk} \cdot \tilde{\sigma}_n([-2^{-k}\eta, 2^{-k}\eta]^d) \geq \eta^d \cdot \phi_\delta(0) > \eta^d/4(2K + 2)^d$$

for every $\delta > 0$ and conclude from (8) that $\tilde{f}$ is recurrent. $\square$